%% file: Distribution_of_Non_Markovian_Coins.tex
\documentclass[11pt,lettersize,reqno,fullpage]{article}
\usepackage{amsmath,amssymb,amsfonts,mathrsfs,verbatim,enumitem,pstricks,amsthm}
\usepackage[all]{xy}
\usepackage{centernot, xr, accents}
\usepackage{hyperref}
 
\include{Macros_for_Probability}

\title{{\LARGE Probability distribution of constrained Random Walks}}
\author{Ritwik Mukherjee }
\date{}
\usepackage{hyperref}
\hypersetup{
	colorlinks,
	citecolor=black,
	filecolor=black,
	linkcolor=blue,
	urlcolor=black
}

\newtheorem{thm}{Theorem}[section]

\newtheorem{rem}[thm]{Remark}

\newtheorem{question}[thm]{Question}

\newtheorem{prpn}[thm]{Proposition}
\newtheorem*{ack}{Acknowledgements}

\begin{document}

\maketitle

\begin{abstract}
In this paper we consider a sequence of $n$ coin tosses, whose outcome depends on the previous   
$n-1$ tosses. In particular, their distribution is not i.i.d. 
We compute the limiting distribution 
of this sequence 
using the method of images.
\end{abstract}

\tableofcontents

\section{Introduction}
We study the following question - suppose a particle is taking a random walk on a 
finite rod. 
What is the chance that it will be at a particular position? 
How much will it deviate on an average? 
We answer this question using the method of images.  
Intuitively, one would expect that in the 
limiting case, the answer should be the fundamental solution of the Heat Equation on a finite rod. 
That is indeed the case. 
The idea we employ is very similar to the 
idea used in \cite{JK}, 
where the authors solve the Schr\"{o}dinger Equation for the particle in a box by explicitly 
computing the propagator.\\ 
\hf \hf To motivate our results, let us first 
recall the  Central Limit Theorem (C.L.T), 
which is for i.i.d random variables. 
\begin{thm} \textbf{(Central Limit Theorem, c.f. \cite{DK})}
\label{clt} 
Let $x$ be any real number and $t$ 
a positive real number.  
Toss a fair coin $n$ times and denote $n_{\mathrm{H}}$ and $n_{\mathrm{T}}$ 
to be the number of heads and tails respectively. Define 
\[ \Sum_{n}^{x,t} := x + \Big( \frac{n_{\mathrm{H}} - n_{\mathrm{T}}}{\sqrt{n}} \Big) \sqrt{t}. \]
Then for any extended real numbers $a$ and $b$ 
\begin{align*}
\lim_{n\lra \infty} \textnormal{Prob}\Big( a\leq  \mathrm{S}_n^{x,t} \leq b \Big) &= \int\limits_{a}^{b} \G(t, x,y) dy,  
\qquad \textnormal{where} \qquad \G(t,x,y) := \frac{1}{\sqrt{2 \pi t}} \exp\Big(-\frac{(y-x)^2}{2 t}\Big). 
\end{align*}
\end{thm}

\section{Main Results}
We will now modify the setup of the C.L.T and consider coin tosses that are not i.i.d. 
The two main theorems we prove in this paper are: 
\begin{thm}
\label{clt_half_line}
Let $x$ and $t$ 
be positive real numbers and 
$n \in \mathbb{N}$ a positive integer.
Toss a coin \textsf{exactly} $n$ times.  
Let $n_{\mathrm{H}}^k$ and $n_{\mathrm{T}}^k$ 
denote the number of heads and tails respectively 
after $k$ tosses. 
Define 
\begin{align}
\Sum_{k}^{x,t} & := x + \Big( \frac{n_{\mathrm{H}}^k - n_{\mathrm{T}}^k}{\sqrt{n}} \Big) \sqrt{t}. 
\label{sum_k_half}
\end{align}
The coin is a ``fair'' coin except for one small difference - for each $k$, the quantity 
$\Sum_k^{x,t}$ has to be strictly greater than zero. 
Then for any extended real numbers $a$ and $b$ such that $(a,b) \subset [0,\infty]$ 
\begin{align*}
\lim_{n\lra \infty} \textnormal{Prob} \Big( a\leq  \mathrm{S}_n^{x,t} \leq b \Big) &= 
\frac{\int\limits_{a}^{b} \H(t, x,y) dy}{\int\limits_{0}^{\infty} \H(t, x,y) dy},  
\qquad \textnormal{where} \qquad \H(t,x,y) := \G(t,x,y) - \G(t,x,-y). 
\end{align*}
\end{thm}

\begin{thm}
\label{clt_box}
Let  $\L$ be a positive real numbers and  $x$  a real number strictly between $0$
and $\L$. 
Let 
$t$ 
be positive real numbers and 
$n \in \mathbb{N}$ a positive integer.
Toss a coin \textsf{exactly} $n$ times. Denote $n_{\mathrm{H}}^k$ and $n_{\mathrm{T}}^k$ 
to be the number of heads and tails respectively after $k$ tosses and define 
\begin{align}
\Sum_{k}^{x,t} & := x + \Big( \frac{n_{\mathrm{H}}^k - n_{\mathrm{T}}^k}{\sqrt{n}} \Big) \sqrt{t}. 
\label{sum_k_box}
\end{align}
The coin is a ``fair'' coin except for one small difference - for each $k$, 
$\Sum_k^{x,t}$ has to be strictly between $0$ and $\L$.
Then for any real numbers $a$ and $b$ such that $(a,b) \subset [0,\L]$
\begin{align*}
\lim_{n\lra \infty} \textnormal{Prob} \Big( a\leq  \mathrm{S}_n^{x,t} \leq b \Big) &= 
\frac{\int\limits_{a}^{b} \K(t, x,y) dy}{\int\limits_{0}^{\L} \K(t, x,y) dy},  
\qquad \textnormal{where} \\ 
\K(t,x,y) & := 
\frac{1}{\mathrm{L}}\sum_{m = -\infty}^{\infty} \sin\Big(\frac{m \pi x }{\mathrm{L}}\Big) 
 \sin\Big(\frac{m\pi y }{\mathrm{L}}\Big) \exp\Big(-\frac{m^2 \pi^2 t}{\mathrm{L}^2}\Big). 
\end{align*}
\end{thm}

\begin{rem}
The reader should observe that  
$\G(t,x,y)$, $\H(t,x,y)$ and $\K(t,x,y)$ are exactly the same as the 
fundamental solution for the Heat equation on an 
infinite rode, semi-infinite rod and finite rod respectively.    
\end{rem}

\begin{rem}
The reader should note that in equations \eqref{sum_k_half} and \eqref{sum_k_box}, the 
denominator is $\sqrt{n}$, not $\sqrt{k}$. The $n$ is chosen before performing the 
sequence of tosses, and then we toss the coin \textsf{exactly} $n$ times. 
For example, let us consider the setup of Theorem \ref{clt_half_line}. Suppose $n=9$, $x=3.01$, $t=1$
and we have tossed a coin three times, getting consecutive tails. On the fourth toss 
we will get a head with probability one; if we got a tail then $\Sum^{x,t}_4$ would be negative.
Secondly, once we have chosen $n=9$, we have decided that we will toss \textsf{exactly} 
$9$ times. 
A similar thing holds for the setup of Theorem \ref{clt_box}, where $\Sum^{x,t}_k$ has to stay 
strictly between $0$ and $\L$. 
\end{rem}

\section{Correspondence between forbidden tosses and allowed tosses}
\ni The basic idea in proving Theorem \ref{clt_half_line} and \ref{clt_box} is that 
we establish a one to one correspondence between forbidden tosses and allowed 
tosses, \textit{after we make a suitable reflection}.  After that, the distribution function 
(in the case of Theorem \ref{clt_box})
appears as an infinite sum. 
In order to evaluate the infinite sum, we write the summand as a Fourier Transform, 
and use the Poisson Summation formula. That produces the Heat Kernel for a finite rod. 
The author learned of this idea from \cite{JK}. \\ 
\hf \hf First, let us set up some notation. We will identify a sequence of 
$n$ tosses with an element 
\[ \omega:= \omega_1 \omega_2 \ldots \omega_n \in \{-1,1\}^n, \]
with $1$ denoting heads and $-1$ denoting tails. 
For each $k\leq n$, define 
\begin{align*}
\Sum_{k}(\omega) &:=  x+ \frac{\sum\limits_{i=1}^{k} \omega_i \sqrt{t}}{\sqrt{n}}.
\end{align*}

\ni Next, given a set $\U \subset \R$ we define 
$\mathcal{P}_n(\U)$ 
to be the set of all $n$-tosses $\omega$ such that 
$\mathrm{S}_{n}(\omega)$ belongs to $\U$, i.e. 
\begin{align}
 \mathcal{P}_n(\U) &:= \{ \omega \in \{-1,1\}^{n}: \mathrm{S}_{n}(\omega) \in 
 \U \}. \label{pn_defn}
\end{align}
Next, given a positive integer $m \leq n$, we define 
the 
following 
two subsets of $\mathcal{P}_n(\U)$:
\begin{align*}
\mathcal{F}^0_{n, m}(\U) &:= \Big\{ \omega \in \mathcal{P}_{n}(\U): \exists i_1 < i_2 <\ldots <i_m 
\in \{0, 1,2, \ldots n\} ,  \\ 
& \qquad \qquad \textnormal{such that} \quad  \Sum_{i_{2k-1}}(\omega) \leq 0, 
~~\Sum_{i_{2k}}(\omega) \geq \L  \qquad \forall ~~ k \in  
\{ 1,2, \ldots, \left \lfloor{\frac{m}{2}}\right \rfloor+1\} \Big\}, \\
%
\mathcal{F}^{\L}_{n,m}(\U) &:= \Big\{ \omega \in \mathcal{P}_n(\U): \exists i_1 < i_2 <\ldots 
<i_m \in \{0, 1,2, \ldots n\} ,  \\ 
& \qquad \qquad \textnormal{such that} \quad  \Sum_{i_{2k-1}}(\omega) \geq \L, 
~~\Sum_{i_{2k}}(\omega) \leq 0  \qquad \forall ~~ 
k \in \{ 1,2, \ldots, \left \lfloor{\frac{m}{2}}\right \rfloor+1\}\Big\}. 
\end{align*}
Observe that 
\begin{align}
\mathcal{F}^0_{n,m}(\U) \cap \mathcal{F}^{\L}_{n,m}(\U) &= 
\mathcal{F}^0_{n,m+1}(\U) \cup \mathcal{F}^{\L}_{n,m+1}(\U) 
\qquad \forall m \in \{1,2 \ldots, n-1\}. \label{bijection_union_intersection}\\
\nonumber 
\end{align}
Finally, given a collection of real numbers $\lambda_1, \lambda_2, \ldots , \lambda_k$, we define  
\begin{align*}
\U_{\lambda_1} &:= \{ x \in \R: 2\lambda_1-x \in \U \} \qquad \textnormal{and} \qquad 
\U_{\lambda_1, \lambda_2, \ldots, \lambda_k} := (\U_{\lambda_1, \lambda_2, \ldots, \lambda_{k-1}})_{\lambda_k}. 
\end{align*}
The set $\U_{\lambda_1}$ is simply the reflection of 
the set $\U$ through the point $x=\lambda_1$.  
Similarly, $\U_{\lambda_1, \lambda_2, \ldots, \lambda_k}$ is obtained by first reflecting 
$\U$ through $x=\lambda_1$, then through $x=\lambda_2$ and so on until $x=\lambda_k$. \\
\hf \hf We are  now ready to 
establish a one to one correspondence 
between forbidden tosses and allowed tosses, after 
we take a suitable reflection. 
\begin{prpn}
\label{forbidden_allowed_bijection}
Let $\U \subset [0, \infty)$ be an open set and $x \in (0, \infty)$. Then
there exists a bijection between $\mathcal{F}^0_{n,1}(\U)$ and 
$\mathcal{P}_n(\U_0)$ if $n$ is sufficiently large.
\end{prpn}

\begin{prpn}
\label{forbidden_allowed_bijection_box}
Let $\U \subset [0, \L]$ be an open set and $x \in (0, \L)$. 
If $n$ is sufficiently large, 
there exists a bijection 
\begin{align}
\mathcal{F}^0_{n,m}(\U) & \iff  \mathcal{P}_n(\U_{0,-\L, -2\L, \ldots, -(m-1)\L }) \qquad 
\textnormal{and}  \label{bj2} \\
\mathcal{F}^{\L}_{n,m}(\U) & \iff  \mathcal{P}_n(\U_{\L,2\L, \ldots, m \L })\label{bj3}
\end{align}
for every $m \in \{1,2, \ldots, n\}$.
\end{prpn}

\textbf{Proof of Proposition \ref{forbidden_allowed_bijection}:} We will explicitly show the bijection. 
Let $\omega \in \mathcal{F}^0_{n,1}(\U)$ and $i_1$ be the smallest integer such that 
$\Sum_{i_1}(\omega) \leq 0$. 
Define the map 
\[ \Phi:\mathcal{F}^0_{n,1}(\U) \lra \mathcal{P}_n(\U_0), \qquad \textnormal{given by} 
\qquad \Phi(\omega) := \omega_1 \ldots \omega_{i_1} \omega_{i_1+1}^* \ldots \omega_{n}^*,\]
where $1^* := -1$ and $-1^* := 1$.
It is easy to see that $\Phi$ maps to $\mathcal{P}_n(\U_0)$ if $n$ is large 
and $\U$ is open. 
Next, suppose $\omega \in \mathcal{P}_{n}(\U_0)$  
and $i_1$ be the smallest integer such that 
$\Sum_{i_1}(\omega) \leq 0$. 
Define the map 
\[ \Psi:\mathcal{P}_n(\U_0) \lra \mathcal{F}^0_{n,1}(\U), \qquad 
\textnormal{given by} 
\qquad \Psi(\omega) := \omega_1 \ldots \omega_{i_1} \omega_{i_1+1}^* \ldots \omega_{n}^*.
\]
It is easy to see that $\Psi$ maps to $\mathcal{F}^0_{n,1}(\U)$ 
if $n$ is sufficiently large and $\U$ is open. 
It is also immediate 
that $\Phi$ and  $\Psi$ are inverses of each other, since  
\[ \Phi \circ \Psi = \textnormal{Id} 
\qquad \textnormal{and} \qquad \Psi \circ \Phi = 
\textnormal{Id}.\]
This proves Proposition \ref{forbidden_allowed_bijection}. \qed \\

\textbf{Proof of Proposition \ref{forbidden_allowed_bijection_box}:} We will prove \eqref{bj2}; 
the proof of \eqref{bj3} is identical. 
Let $\omega \in \mathcal{F}^0_{n,m}(\U)$ and 
$i_1, i_2, i_3, \ldots, i_m$ be the unique collection of $m$ integers 
such that  
$i_1$ is the smallest integer such that 
$\Sum_{i_1}(\omega) \leq 0$ 
and 
for all $r \in \{2, \ldots, m\}$, $i_r$ is the smallest integer such that  
$i_r > i_{r-1}$ and $\Sum_{i_{r}}(\omega)$ is either less than 
or equal to zero or greater than or equal to $\L$, depending on whether $r$ is odd 
or even respectively. 
Now we define the map 
\begin{align*}
\Phi& :\mathcal{F}^0_{n,m}(\U) \lra \mathcal{P}_{n}(\U_{0,-\L, -2\L, \ldots, -(m-1)\L }), 
\qquad \textnormal{given by} \\
\Phi(\omega) & := \omega_1 \ldots \omega_{i_1} \omega_{i_1+1}^* \ldots \omega_{i_2}^* 
\omega_{i_2+1} \ldots \omega_{i_3} \ldots.  
\end{align*}
It is easy to see that $\Phi$ maps to $\mathcal{P}_{n}(\U_{0,-\L, -2\L, \ldots, -(m-1)\L })$  
if $n$ is sufficiently large and $\U$ is open.\\
\hf \hf Next, suppose $\omega \in \mathcal{P}_{n}(\U_{0,-\L, -2\L, \ldots, -(m-1)\L })$. 
Let $i_1, i_2, i_3, \ldots, i_m$ be the unique collection of $m$ integers 
such that  
$i_1$ is the smallest integer such that 
$\Sum_{i_1}(\omega)\leq 0$ 
and 
for all $r \in \{2, \ldots, m\}$, $i_r$ is the smallest integer such that  
$i_r > i_{r-1}$ and $\Sum_{i_{r}}(\omega)$ is less than 
or equal to $ -(r-1)\L$.  
Now we define the map 
\begin{align*}
\Psi &:\mathcal{P}_n(\U_{0,-\L, -2\L, \ldots, -(m-1)\L }) \lra 
\mathcal{F}^0_{n,m}(\U), 
\qquad \textnormal{given by} \\
\Psi(\omega) & := \omega_1 \ldots \omega_{i_1} \omega_{i_1+1}^* \ldots \omega_{i_2}^* 
\omega_{i_2+1} \ldots \omega_{i_3} \ldots.  
\end{align*}
It is easy to see that $\Psi$ maps to $\mathcal{F}^0_{n,m}(\U)$ 
if $n$ is sufficiently large and $\U$ is open. 
It is also immediate 
that $\Phi$ and $\Psi$ are inverses of each other, since 
\[ \Phi \circ \Psi = \textnormal{Id}
\qquad \textnormal{and} \qquad \Psi \circ \Phi = 
\textnormal{Id}.\]
This proves \eqref{bj2}. The proof of \eqref{bj3} is identical.  \qed 

\section{Proofs of the main results}
\ni We are now ready to prove Theorem \ref{clt_half_line} and \ref{clt_box}. \\

\textbf{Proof of Theorem \ref{clt_half_line}:} First, we observe that for sufficiently large $n$
\begin{align*}
\textnormal{Prob}\Big( a\leq  \mathrm{S}_n^{x,t} \leq b \Big) &= 
\frac{|\mathcal{P}_n((a,b))|- |\mathcal{F}^0_{n,1}((a,b))|}{| \mathcal{P}_{n}((0,\infty))|
-|\mathcal{F}^0_{n,1}((0,\infty))|}  \\
& = \frac{|\mathcal{P}_{n}((a,b))|- |\mathcal{P}_n((-b,-a))|}{| \mathcal{P}_n((0,\infty))| 
-|\mathcal{P}_n((-\infty,0))| } \qquad \textnormal{by 
Proposition \ref{forbidden_allowed_bijection}} \\ 
& = \frac{(|\mathcal{P}_n((a,b))|- |\mathcal{P}_n((-b,-a))|)\times 2^{-n}}{(|\mathcal{P}_n((0,\infty))| 
-|\mathcal{P}_n((-\infty,0))|)\times 2^{-n}} \\
\implies \lim_{n\lra \infty} \textnormal{Prob}\Big( a\leq  \mathrm{S}_n^{x,t} \leq b \Big) &= 
\frac{\int\limits_{a}^b \G(t,x,y) dy - \int\limits_{-b}^{-a} \G(t,x,y) dy}{\int\limits_{0}^{\infty} \G(t,x,y) dy - 
\int\limits_{-\infty}^{0} \G(t,x,y)dy} \qquad \textnormal{by C.L.T} \\
& = \frac{\int\limits_{a}^{b} \H(t,x,y) dy}{\int\limits_{0}^{\infty} \H(t,x,y) dy}. \qquad \textnormal{\qed}
\end{align*}

\textbf{Proof of Theorem \ref{clt_box}:} 
Let $\mathcal{F}_n((a,b))$ denote the set of all forbidden $n$-tosses in the setup of this Theorem. 
Note that 
\begin{align}
\mathcal{F}_n((a,b)) &= \mathcal{F}^0_{n,1}((a,b)) 
\cup \mathcal{F}^{\L}_{n,1}((a,b)).  \label{forbidden_f0_f1_union} 
\end{align}
Using equations \eqref{forbidden_f0_f1_union}, \eqref{bijection_union_intersection}, Proposition 
\ref{forbidden_allowed_bijection_box} and the Inclusion Exclusion Principle, we conclude that  
\begin{align*}
|\mathcal{P}_{n}((a,b))|- |\mathcal{F}_{n}((a,b))| &= \sum_{m=-\infty}^{\infty} 
\mathcal{P}_n((2m\L +a, 2m\L +b)) -\mathcal{P}_n((2m\L-b, 2m\L-a)).  
\end{align*}
Note that the above expression is actually a finite sum; when $m$ is sufficiently large the terms will 
become zero. Hence, we conclude that  for large $n$,
\begin{eqnarray}
\textnormal{Prob}\Big( a\leq  \mathrm{S}_n^{x,t} \leq b\Big) = 
\frac{|\mathcal{P}_{n}((a,b))|- |\mathcal{F}_n((a,b))|}{|\mathcal{P}_n((0,\L))|- |\mathcal{F}_n((0,\L))|} 
\nonumber \\ 
 = \frac{\sum\limits_{m=-\infty}^{\infty} \mathcal{P}_n((2m\L+a, 2m\L+b)) - 
\mathcal{P}_n((2m\L -b, 2m\L -a)) }{\sum\limits_{m=-\infty}^{\infty} \mathcal{P}_n((2m\L, 2m\L+\L)) - 
\mathcal{P}_n((2m\L -\L, 2m\L))} \nonumber\\
 = \frac{\sum\limits_{m=-\infty}^{\infty} \Big( \mathcal{P}_n((2m\L+a, 2m\L+b)) - 
\mathcal{P}_n((2m\L -b, 2m\L -a)) \Big) \times 2^{-n} }
{\sum\limits_{m=-\infty}^{\infty} \Big( \mathcal{P}_n((2m\L, 2m\L+\L)) - 
\mathcal{P}_n((2m\L -\L, 2m\L)) \Big) \times 2^{-n}} \label{dct1}\\
\implies \lim_{n\lra \infty} \textnormal{Prob}\Big( a\leq  \mathrm{S}_n^{x,t} \leq b \Big) =  
\frac{\int\limits_{a}^{b} \mathcal{K}(t,x,y)dy}{\int\limits_{0}^{\L} \mathcal{K}(t,x,y)dy}, \qquad \textnormal{where} 
\label{dct2} \\
\mathcal{K}(t,x,y)  := \sum\limits_{m=-\infty}^{\infty} \mathcal{G}(t,x,y+2m\L) - \mathcal{G}(t,x,-y+2m\L). \nonumber 
\end{eqnarray}
To go from \eqref{dct1} to \eqref{dct2},   
we interchanged the order of limit and summation; 
this is justified in section \ref{dct_justify_appendix}. 
Finally, we evaluate $\K$ by writing $\G$ as a Fourier Transform and 
using the Poisson Summation Formula. Let  
\begin{align*}
f(k) &:= \exp(-2\pi^2tk^2 + 2\pi k\i (x-y))-\exp(-2\pi^2 t k^2 + 2\pi k\i (x+y)). 
\end{align*}
Then we get that 
\begin{align*}
\K(t,x,y) &= \sum\limits_{m= -\infty}^{\infty} \int\limits_{-\infty}^{\infty} f(k) \exp(-4\pi \i m k \L) dk \\ 
                          &= \frac{1}{2 \L} \sum_{m = -\infty}^{\infty} f\Big(\frac{m}{2 \L}\Big) \qquad 
                          \textnormal{(by the Poisson
                          Summation Formula, c.f. \cite{ES})} \\ 
                          &= \frac{1}{\mathrm{L}}\sum_{m = -\infty}^{\infty} \sin\Big(\frac{m \pi x }{\mathrm{L}}\Big) 
 \sin\Big(\frac{m \pi y }{\mathrm{L}}\Big) \exp\Big(-\frac{m^2 \pi^2 t}{\mathrm{L}^2}\Big). \qquad \textnormal{\qed}
\end{align*}



\section{Justifying the interchange of limit and summation}
\label{dct_justify_appendix}
\ni We now justify the interchange of limit and summation 
that was used in going from \eqref{dct1} to \eqref{dct2}. 
First, we state a couple of well known inequalities. 
\begin{thm}
\label{Stirling}
(Stirling's Inequality)
For all $n \in \mathbb{N}$, we have  
\begin{align*}
\sqrt{2 \pi} n^{n+ \frac{1}{2}} e^{-n} & \leq n! \leq e n^{n+ \frac{1}{2}} e^{-n}.  \label{stirling_ineq}
\end{align*}
\end{thm}


\begin{thm}
\label{Bernoulli}
(Bernoulli's Inequality) 
Let $x,r \in \mathbb{R}$, such that $x >-1 $ and $r \geq 1$. Then 
\[ \frac{1}{(1+x)^r} \leq \frac{1}{1+r x}. \]
\end{thm}

\begin{prpn}
Let $\U$ be a subset of $\R$ and 
$\mathcal{P}_n(\U)$ be as defined in  \eqref{pn_defn}.  
Define
\[ \alpha_{m,n} := \frac{|\mathcal{P}_n((2m\L+ a, 2m\L+b))|}{2^n} \qquad \textnormal{and} 
\qquad 
\beta_{m,n} := \frac{|\mathcal{P}_n((2m\L -b, 2m\L-a))|}{2^n}. \]
Then there exists a constant $C$ (independent of $m$ and $n$), 
such that 
\begin{align}
\alpha_{m,n} & \leq  \frac{C}{\Big(1- \frac{(2m\L +a-x)^2}{2t} \Big)^2 
\Big(1-\frac{(2m\L +a-x)^2}{t} \Big)} \label{alpha_mn} \\
\beta_{m,n} & \leq \frac{C}{\Big(1- \frac{(2m\L -b-x)^2}{2t} \Big)^2 
\Big(1-\frac{(2m\L -b-x)^2}{t} \Big)} \label{beta_mn}
\end{align}
Furthermore, 
\begin{align}
\lim_{n\lra \infty} \sum_{m=-\infty}^{\infty} (\alpha_{m,n} - \beta_{m,n})  &= 
\sum_{m=-\infty}^{\infty} \lim_{n\lra \infty}(\alpha_{m,n} - \beta_{m,n}). \label{dct_justify_again}
\end{align}
\end{prpn}
\textbf{Proof:} We note that 
\begin{align*}
\alpha_{m,n} &\leq  \sum_{i=0}^{b-a+1}\binom{n}{ \left\lfloor{\frac{n}{2} - 
\Big( \frac{2m\L +a+i-x}{2 \sqrt{t}} } \Big)\sqrt{n} \right\rfloor} \\
& \leq (b-a+1) \binom{n}{ \left\lfloor{\frac{n}{2} - 
\Big( \frac{2m\L +a-x}{2 \sqrt{t}} } \Big)\sqrt{n} \right\rfloor} \\
& \leq \frac{C^{\prime}}{\sqrt{n} \Big(1- \frac{(2m\L +a-x)^2}{t} \frac{1}{n} \Big)^{\frac{n}{2} + \frac{1}{2}} 
\Big(1- \frac{(2m\L +a-x)^2}{\sqrt{t}} \frac{1}{\sqrt{n}} \Big)^{-\Big(\frac{(2m\L +a-x)^2}{2 \sqrt{t}}\Big) \sqrt{n}}
\Big(1+\frac{(2m\L +a-x)^2}{\sqrt{t}} \frac{1}{\sqrt{n}} \Big)^{\Big(\frac{(2m\L +a-x)^2}{2 \sqrt{t}}\Big) \sqrt{n}}
}\\
& \textnormal{(using Theorem \ref{Stirling})} \\
& \leq \frac{C}{\Big(1- \frac{(2m\L +a-x)^2}{2t} \Big)^2 
\Big(1-\frac{(2m\L +a-x)^2}{t} \Big)}  \qquad \textnormal{(using Theorem \ref{Bernoulli}),}
\end{align*}
where $C^{\prime}$ is some constant. 
Inequality \eqref{beta_mn} follows from \eqref{alpha_mn} by replacing $a$ by $-b$ and $b$ by 
$-a$. Finally, \eqref{dct_justify_again} follows from \eqref{alpha_mn}, \eqref{beta_mn}, the 
triangle inequality and the Dominated Convergence Theorem (since the rhs of 
\eqref{alpha_mn} and \eqref{beta_mn} are infinite summable).  
Hence, \eqref{dct1} implies \eqref{dct2}. \qed 

\section{Further remarks}
\label{future}

\ni We end this paper with a few natural questions that one can investigate in the future. 

\begin{question}
Can this idea of ``method of images'' be used to find the limiting distribution of constrained 
random 
walks defined on a free group (or a finitely generated group)? 
After all, it seems quite natural to extend the idea of taking a suitable reflection on any group.  
\end{question}

\begin{question}
One can consider higher dimensional analogues of Theorem \ref{clt_half_line} and \ref{clt_box}. 
Is the distribution function the same as the fundamental solution of the corresponding higher dimensional 
Heat Equation? 
\end{question}

\begin{question}
Do the random variables considered in Theorem \ref{clt_half_line} and \ref{clt_box} satisfy a 
Large Deviation Principle? If yes, is it possible to compute the rate function explicitly? 
\end{question}

\begin{ack}
The author is grateful to Amritanshu Prasad and Vamsi Pingali for several useful discussions. 
\end{ack}


\end{document}

%% file: Macros_for_Probability.tex
\def \ni{\noindent}
\def \hf{\hspace*{0.5cm}}                      

\def\lra{\longrightarrow}

\def\bge{\begin{equation}}                
\def\ede{\end{equation}}                
\def\bgd{\begin{displaymath}}         
\def\edd{\end{displaymath}}            
\def\bgee{\begin{equation*}}           
\def\edee{\end{equation*}}           

\def \R{\mathbb{R}}
\def \L{\mathrm{L}}

\def \G{\mathcal{G}}
\def \H{\mathcal{H}} 
\def \K{\mathcal{K}} 

\def \Sum{\mathrm{S}}

\def \U{\mathcal{U}}
\def \i{\sqrt{-1}}

%% file: Distribution_of_Non_Markovian_Coins.bbl
\begin{thebibliography}{99}

\bibitem{JK} W.Janke, H.Kleinert,   ~\textit{Summing Paths for Particle in a Box}, \\
Lettere Al Nuovo Cimento, Vol. $25$, N. $10$.    

\bibitem{DK} D.Khoshnevisan, ~\textit{Probability}, 
Graduate Studies in Mathematics, Volume $80$. 

\bibitem{ES} E.M.Stein, R.Shakarchi, ~\textit{Fourier Analysis}, 
Princeton Lectures in Analysis.



\end{thebibliography}
